\documentclass[12pt, twoside, a4paper]{article}
\usepackage[T1]{fontenc}
\usepackage{makeidx}\usepackage{amssymb}
\usepackage{amsmath,amssymb,amsopn,amsfonts,mathrsfs,amsbsy,amscd}
\usepackage[usenames,dvipsnames]{pstricks} \usepackage{epsfig}

\newcommand{\too}{\longrightarrow}
\newcommand{\om}{\omega}
\newcommand{\esp}{\quad\mbox{and}\quad}

\newcommand{\di}{\displaystyle}
\newcommand{\Om}{\Omega}
\newcommand{\na}{\nabla}

\newcommand{\al}{\alpha}

\newcommand{\Ga}{\Gamma}
\newcommand{\e}{\epsilon}

\newcommand{\la}{\lambda}

\font\bb=msbm10

\def\R{\hbox{\bb R}}

\newtheorem{theorem}{Theorem}[section]

\newtheorem{e-proposition}[theorem]{Proposition}

\newtheorem{e-definition}[theorem]{Definition\rm}

\newtheorem{example}{\it Example\/}

\title{Symplectic structures on the  tangent bundle of a smooth manifold
}
\author{Abdelhaq Abouqateb, Mohamed Boucetta, Aziz Ikemakhen}\date{}

\begin{document}\maketitle

\begin{abstract}
 We give a method to lift $(2,0)$-tensors fields on a manifold $M$ to build symplectic
forms on $TM$. Conversely,  we show that any  
symplectic
form $\Om$ on
$TM$   is symplectomorphic, in a neighborhood of the zero section, to a symplectic form
built naturally from three $(2,0)$-tensor fields   associated to $\Om$.

\end{abstract}

\section{Introduction}
\label{section1}

The geometry of the tangent bundle $TM$ of a Riemannian manifold
 with the Sasaki metric has been extensively studied since the 60's (see
\cite{boris,sasaki}).
To
overcome the rigidity of these metrics many others metrics generalizing Sasaki metrics
where introduced and studied (see \cite{kowal}).  Given a Riemannian manifold
$(M,g)$, the basic idea  behind the
construction of Riemannian metrics on $TM$ from the Riemannian metric on $M$ is
to
use the Levi-Civita connection of $g$ to get a splitting of $TTM=\mathcal{V}M\oplus
\mathcal{H}M$  and to lift the metric $g$ to $TM$ by the mean of this splitting. It is
natural, in order to construct symplectic structures on the tangent bundle, to use the
same approach in the case where Riemannian metrics are replaced by differential
2-forms or, more generally, $(2,0)$-tensor fields. This point of view has been adopted in 
\cite{janyska0,janyska,kurek} motivated by the classification of  "natural" symplectic
forms on the tangent bundle.  In this paper, we address the two following situations:
\begin{enumerate}
 \item Starting with a manifold $M$ endowed with two differential 2-forms $\om_0,\om_1$,
 a $(2,0)$ tensor field $A$
 and a linear connection $\na$, we construct a natural  differential 2-form $\Om$
on $TM$ involving  the $\om_i$, $A$ and the splitting of $TTM$ induced by $\na$. We give
then
the sufficient and necessary conditions on $\om_i$, $A$ and $\na$ for which $\Om$ is
symplectic
(see Proposition \ref{proposition1}). Among these conditions, $(A,\na)$ satisfy an
equation which is
known as Codazzi equation (see \cite{delanoe,shima}) when $A$ is a Riemannian metric and
$\na$ flat. We give in
Proposition \ref{propisition2} some equivalent assertions to this equation.
\item Conversely, to any symplectic form $\Om$ on $TM$,  we
associate two differential 2-forms $\om_{11},\om_{22}\in\Om^2(M)$ and  
$(2,0)$-tensor fields $A$ on $M$.  We show that, for any choice of a connection $\na$ on
$M$, $\Om$   is symplectomorphic, near the zero section,  to a symplectic form built
from $(\na,\om_{11},\om_{22},A)$ in a  way described in (i)
(see Theorem \ref{main}). 
\end{enumerate}The paper is organized as follows. In Section \ref{sections} we state our
main result and
in Section \ref{section2}, we
define the
lift of $(2,0)$-tensor fields to the tangent bundle by the mean of a linear connection
and we prove Propositions \ref{propisition2}-\ref{proposition1}. 
Section \ref{section3} is devoted to a proof of Theorem \ref{main} which is mainly based
on a version of the classical Darboux's Theorem. 

\section{Statement of the main result}\label{sections}

Let $M$ be a	manifold and $\pi:TM\too M$ its tangent bundle. We denote by $\imath:
M\too TM$ the zero section and by $\mathcal{V}M=\ker d\pi$ the vertical subbundle of
$TTM$. For any  $x\in M$ and $u\in T_xM$ there is a natural isomorphism
$\tau_{(x,u)}:T_xM\too \mathcal{V}_uM$ given by
$\tau_{(x,u)}(v)=\frac{d}{dt}_{|t=0}(u+tv).$
For any vector field $X$ on $M$,
we define its {\it vertical lift} $X^{v}$ which is the vector field on $TM$ given by
$X^{v}(x,u)=\tau_{(x,u)}(X_x).$
On the other hand, for any $x\in M$, we have
$T_{0_x}TM=\mathcal{V}_{0_x}M\oplus
\imath_*(T_xM),$ where $0_x$ is the null vector of $T_xM$. For any
differential 2-form $\Om$ on $TM$, we
associate three $(2,0)$-tensor fields   $\om_{11},\om_{22},A$ on $M$  by putting
\begin{eqnarray}
 \om_{11}&=&\imath^*\Om,\label{eqse13},\;
A(u,v)(x)=\Om(\tau_{(x,0_x)}(u),\imath_*(v))_{0_x},\\
\om_{22}(u,v)(x)&=&\Om(\tau_{(x,0_x)}(u),\tau_{(x,0_x)}(v))_{0_x}.\nonumber
\end{eqnarray}
Suppose now that $M$ carries a linear connection $\na$. This define an horizontal
distribution on $TM$ as follows:
$$\mathcal{H}_{(x,u)}M=\left\{\frac{d}{dt}_{|t=0}P_{c_x}^u(t),c_x\in
C^\infty(]-\e,\e[,M)\;\mbox{and}\; c_x(0)=x\right\}$$where $P_x^u(t):]-\e,\e[\too TM$ is
the parallel transport with respect to $\na$ of $u$ along the curve $c_x$. The linear map
$T_{(x,u)}\pi:\mathcal{H}_{(x,u)}M\too T_xM$ is an isomorphism and hence
$TTM=\mathcal{V}M\oplus \mathcal{H}M.$
For any vector field $X$ on $M$, we define its horizontal lift $X^h$ by
$X^h_{(x,u)}=\left(T_{(x,u)}\pi\right)^{-1}(X_x).$
For any $\om\in\Om^2(M)$ we define the differential 1-form $\la^{\om,\na}$ on $TM$ by
putting,
for any vector field $X$ on $M$ and for any $u\in TM$,
\begin{equation}
 \la^{\om,\na}(X^h)=0\esp
\la^{\om,\na}(X^v)(u)=\frac12\om(u,X).\label{eqse16}
\end{equation}

We can now state our main result.
\begin{theorem}\label{main}   Let $\Om$ be a
symplectic form on $TM$. Let $\om_{11},A,\om_{22}$
 be the associated $(2,0)$-tensor fields given by (\ref{eqse13}).
 Then, for any linear connection $\na$,
 there exists two open  neighborhoods $\mathcal{N}_1$ and $\mathcal{N}_2$ of the zero
section in $TM$
and a diffeomorphism
$\phi:\mathcal{N}_1\too\mathcal{N}_2$ such that
$$\phi_{|\imath(M)}=\mathrm{Id}_{\imath(M)} \esp
\phi^*\Om=\pi^*\om_{11}+(A^\flat)^*(d\la)+d\la^{\om_{22},\na},$$
where $\la$ is the Liouville 1-form on $T^*M$ and
$A^\flat:TM\too T^*M$ is given
by $A^\flat(u)=A(u,.)$.

\end{theorem}

\section{Lift of $(2,0)$-tensor fields on $M$ to  symplectic forms on
$TM$}\label{section2}
Let $M$ be a manifold endowed with a linear connection $\na$
 and
$TTM=\mathcal{V}M\oplus
\mathcal{H}M$ the associated splitting. 
Let $(x^1,\ldots,x^n)$ be a local coordinates system on $M$ and
$(x^1,\ldots,x^n,u^1,\ldots,u^n)$ the corresponding coordinates system on $TM$. Let
$(\Ga_{ij}^k)$ the Christoffel's symbols of $\na$ defined by
$\na_{\partial_{x^i}}\partial_{x^j}=\sum_{k=1}^n\Ga_{ij}^k\partial_{x^k}.$
If $X=\sum_{i=1}^n X^i\partial_{x^i}$  then
\begin{eqnarray}
 X^h&=&\sum_{i=1}^n X^i\partial_{x^i}-\sum_{i,j,k}\Ga_{ij}^ku^iX^j\partial_{u^k}\esp
X^v=\sum_{i=1}^n X^i\partial_{u^i}.\label{eqse21}
\end{eqnarray}
We deduce easily from these formulas:
\begin{eqnarray}\label{eqse22}
 \;[X^h,Y^h]&=&[X,Y]^h-\left(R(X,Y)u)\right)^v,\quad[X^v,Y^v]=0,\\
\;[X^h,Y^v]&=&(\na_XY)^v,\nonumber
\end{eqnarray}
where $R$ is the curvature of $\na$ given by
$R(X,Y)=\na_{[X,Y]}-\left(\na_X\na_Y-\na_Y\na_X\right).$
Let $\om_0,\om_1$ two differential 2-forms on $M$ and $A$ a $(2,0)$-tensor field.
We define $\Om\in\Om^2(TM)$ by
\begin{eqnarray}\Om(X^v,Y^v)&=&\om_0(X,Y)\circ\pi,\;
\Om(X^h,Y^h)=\om_1(X,Y)\circ\pi,\;
\Om(X^v,Y^h)=A(X,Y)\circ\pi,\nonumber\\
\Om(X^h,Y^v)&=&-A(Y,X)\circ\pi.\label{eqse23}\end{eqnarray}
We call $\Om$ the {\it lift} of $(\na,\om_0,\om_1,A)$. This notions appeared in
\cite{janyska0} when $\om_0=\om_1=0$ and $A$ is a Riemannian metric.
It is obvious that $\Om$ is nondegenerate iff, for any local
coordinates system $(x^1,\ldots,x^n)$ the matrix 
$\left(\begin{array}{cc}P&M\\-M^t&Q\end{array}\right),$ where
$P=\left(\om_{11}(\partial_{x^i},\partial_{x^j})\right)_{1\leq i,j\leq n},\;
Q=\left(\om_{22}(\partial_{x^i},\partial_{x^j})\right)_{1\leq i,j\leq n},\;
M=\left(A(\partial_{x^i},\partial_{x^j})\right)_{1\leq i,j\leq n}$ is invertible. A
direct computation using (\ref{eqse22}) gives:
\begin{eqnarray}
 d\Om(X^h,Y^h,Z^h)&=&d\om_1(X,Y,Z)\circ\pi+\oint
A(R(X,Y)u,Z)\circ\pi,\nonumber\\
d\Om(X^v,Y^v,Z^v)&=&0,\nonumber\\
d\Om(X^v,Y^v,Z^h)&=&\na_{Z}\om_0(X,Y)\circ\pi,\label{closed}\\
d\Om(X^h,Y^h,Z^v)
&=&-\na_{X}A(Z,Y)\circ\pi+\na_{Y}A(Z,X)\circ\pi+A(Z,\tau(X,Y))\nonumber\\&&+\om_0(R(X,Y)u,
Z)\circ\pi,\nonumber
\end{eqnarray}where $\tau$ is the torsion of $\na$ given by
$\tau(X,Y)=[X,Y]-\na_XY+\na_YX.$
We call the equation \begin{equation}\label{codazzi}
 \na_{X}A(Z,Y)-\na_{Y}A(Z,X)=A(Z,\tau(X,Y))
\end{equation}
 Codazzi equation. Indeed, when $\na$ is torsion free
and $A$ is a pseudo-Riemannian metric, we recover the Codazzi equation known in the
context of Hessian manifolds (\cite{shima}). It appeared also in
 \cite{delanoe}. 
The following result is a generalization both of a result of Delano\"e \cite{delanoe} and
a result by Janyska in \cite{janyska0,janyska}.
\begin{e-proposition}\label{propisition2} Let $(M,\na)$ be a manifold endowed with a 
connection and $A$ a nondegenerate $(2,0)$-tensor field. Let $\Om$ be the lift of
$(\na,0,0,A)$. Then the following assertions
are equivalent:
\begin{enumerate}
 \item $(A,\na)$ satisfies Codazzi equation (\ref{codazzi}).
\item $\Om=(A^{\flat})^*(d\la),$ where $\la$ is the Liouville 1-form on $T^*M$
and $
A^\flat:TM\too T^*M$ is given by $A^\flat(u)=A(u,.)$.  
\item $\Om$ is symplectic.
\item $A^\flat(\mathcal{H}M)$ is Lagrangian with respect to $d\la$, where $\la$ is the
Liouville 1-form on $T^*M$ and $\mathcal{H}M$ is the horizontal distribution associated
to $\na$.

\end{enumerate}

\end{e-proposition}

{\bf Proof.} Remark first that since $A$ is nondegenerate then $\Om$ is nondegenerate. We
choose a local coordinates system $(x^i)_{i=1}^n$ and we denote by
$(x^i,u^i)_{i=1}^n$ and $(x^i,p^i)_{i=1}^n$ the corresponding coordinates on $TM$ and
$T^*M$ respectively. The Liouville 1-form is given by 
$\la=\sum_{i=1}^np^idx^i$and $A^\flat$ is given by 
$A^\flat(x^1,\ldots,x^n,u^1,\ldots,u^n)=(x^1,\ldots,x^n,P^1,\ldots,P^n),$ where
$\di P^i=\sum_{j=1}^nu^jA_{ji}$, {with}  $A_{ij}=A(\partial_{x^i},\partial_{x^j}).$
Thus
$(A^\flat)^*(d\la)=\sum_{i=1}^ndP^i\wedge dx^i.$
By using (\ref{eqse21}), we get for $i=1,\ldots,n$,
\begin{eqnarray*}
 \partial_{x^l}^h&=&\partial_{x^l}-\sum_{j,k}\Ga_{jl}^ku^j\partial_{u^k},
\esp
\partial_{x^l}^v=\partial_{u^l}.
\end{eqnarray*}
Thus
$\partial_{x^l}^h(P^s)=\sum_{j=1}^nu^j\left(\partial_{x^l}(A_{js})
-\sum_{k}\Ga_{jl}^kA_{ks}\right).$ Now 
\begin{eqnarray*}
(A^\flat)^*(d\la)\left(\partial_{x^l}^v,\partial_{x^s}^v\right)&=&0,\\
(A^\flat)^*(d\la)\left(\partial_{x^l}^h,\partial_{x^s}^h\right)&=&
\sum_{i=1}^n\left(\partial_{x^l}^h(P^i)\partial_{x^s}^h(x^i)
-\partial_{x^l}^h(x^i)\partial_{x^s}^h(P^i)\right)\\
&=&\partial_{x^l}^h(P^s)-\partial_{x^s}^h(P^l)\\
&=&\sum_{j=1}^nu^j\left(\partial_{x^l}(A_{js})-\partial_{x^s}(A_{jl})
-\sum_{k}\left(\Ga_{jl}^kA_{ks}-\Ga_{js}^kA_{kl}\right)\right)\\
(A^\flat)^*(d\la)\left(\partial_{x^l}^v,\partial_{x^s}^h\right)&=&
\sum_{i=1}^n\left(\partial_{x^l}^v(P^i)\partial_{x^s}^h(x^i)
-\partial_{x^v}^h(x^i)\partial_{x^s}^h(P^i)\right)
=\partial_{x^l}^v(P^s)=
A_{ls}.
\end{eqnarray*}
This shows that $(i)$, $(ii)$ and $(iv)$ are equivalent. Moreover, $(ii)$ implies $(iii)$
obviously and the expression $d\Om(X^h,Y^h,Z^v)$ given in (\ref{closed}) shows that 
$(iii)$ implies $(i)$ . 
\hfill$\square$

\begin{e-proposition}\label{proposition1} The differential 2-form $\Om$ is closed if
and
only if the
following relations hold:
\begin{enumerate}
 \item $d\om_1=0$,
 $\na\om_0=0$ and, for any $X,Y,Z,T$,  $\om_0(R(X,Y)Z,T)=0$,
\item $(\na,A)$ satisfy the Codazzi equation (\ref{codazzi}).

\end{enumerate}

\end{e-proposition}
{\bf Proof.} If $\Om$ is closed then, according to the relations (\ref{closed}), $(i)$
and $(ii)$ hold. Conversely, write $\Om=\Om_1+\Om_2$ where $\Om_1$ is the lift of
$(\na,\om_0,\om_1,0)$ and $\Om_2$ is the lift of $(\na,0,0,A)$. So if $(i)$ and $(ii)$
hold then, by Proposition \ref{propisition2}, $\Om_2=(A^\flat)^*(d\la)$ which is closed.
Hence $\Om$ is closed iff $\Om_1$ is closed which is guaranteed by $(i)$.\hfill$\square$

Let us give some situations where we can use Proposition \ref{proposition1} or
Proposition \ref{propisition2}  to build
symplectic forms on the tangent bundle or Lagrangian horizontal distribution on $T^*M$.
\begin{example}
\begin{enumerate}
 \item Let $(M,\om_1)$ be a symplectic manifold, $\om_0$ a nondegenerate 2-form on $M$
  and $\na$ a flat connection such that $\na\om_0=0$. According to
Proposition \ref{proposition1} the lift of  $(\na,\om_{0},\om_{1},0)$
 is a symplectic form on $TM$. 
\item  Let $G$ be  a Lie
group, $\om_1$ a left invariant symplectic form on $G$, $\om_0$ a nondegenerate right
invariant 2-form on $G$,  $\na$ the flat connection on $G$ satisfying $\na X=0$ for any
right invariant vector field.
Then, according  to
Proposition \ref{proposition1}, the lift of $(\na,\om_{0},\om_{1},0)$ is a symplectic form
on $TG$.
\item Let $M$ be a manifold, $\na$ a  connection on $M$ and $\al$ a differential
1-form on $M$.
Put
$A(X,Y)=\na_Y\al(X).$ One can check easily that
$$ \na_{X}A(Z,Y)-\na_{Y}A(Z,X)=A(Z,\tau(X,Y))+\al(R(X,Y)Z).$$
So if $\na$ is flat then $(A,\na)$ satisfy Codazzi equation. By  choosing $\al$ such that
for any coordinates system $(x^1,\ldots,x^n)$ on $M$ the matrix
$\left(\na_{\partial_{x^i}}\al(\partial_{x^j})\right)_{1\leq i,j\leq n}$ is invertible and
by using   
Proposition \ref{propisition2}, we get that $A^\flat(\mathcal{H})$ is a Lagrangian  
distribution with respect to $d\la$.\\
We give now an example of $\na$ and $\al$ satisfying the conditions above. We consider
$\R^n$ with its canonical $\na$ and let $B=(b_{ij})$
be an invertible $n$-matrix. Put, for $i=1,\ldots,n$,
$\al(\partial_{x^i})=\exp\left(\sum_{k=1}^nb_{ki}x^k\right).$
An easy computation gives that 
$$\left(\na_{\partial_{x^i}}\al(\partial_{x^j})\right)_{1\leq i,j\leq n}=BD$$where $D$ is
the diagonal matrix with entries $\al(\partial_{x^i})$, $i=1,\ldots,n$.
\item Let $(M,g)$ be a pseudo-Riemannian manifold and $\na$ the Levi-Civita connection of
 $g$. According to Proposition \ref{propisition2},  the lift of
$(\na,0,0,g)$  is $(g^\flat)^*(d\la)$. This situation was pointed out in
\cite{janyska0,janyska}. Moreover, $g^\flat(\mathcal{H}M)$ is an horizontal
Lagrangian distribution.
\item Let $(M,\om)$ be a symplectic manifold and $\na$ a torsion free connection such
that $\na\om=0$. It is a well-known result that there are  many such  connections (see
\cite{cahen}). According to Proposition \ref{propisition2},  the lift of
$(\na,0,0,\om)$  is $(\om^\flat)^*(d\la)$ and $\om^\flat(\mathcal{H}M)$ is an horizontal
Lagrangian distribution.

\end{enumerate}

\end{example}

\section{Proof of Theorem \ref{main}}\label{section3}
The proof of Theorem \ref{main} is based on the  following version of the classical
Darboux's theorem (see \cite{macduff}).
\begin{theorem}\label{darboux}
 Let $V$ be a  smooth manifold and $\om_1,\om_2\in\Om^2(M)$ are closed. Suppose
that $N$ is
a submanifold of $V$ such that for any $q\in N$, $\om_1(q)=\om_2(q)$ and
$\om_1(q),\om_2(q)$ are non-degenerate. Then there exists two open neighborhoods
$\mathcal{N}_1$ and $\mathcal{N}_2$ of $N$ and a diffeomorphism
$\phi:\mathcal{N}_1\too\mathcal{N}_2$ such that
$$\phi_{|N}=\mathrm{Id}_N\esp \phi^*\om_2=\om_1.$$
\end{theorem}

{\bf Proof of Theorem \ref{main}.} Put $V=TM$, $N=\imath(M)$, $\om_1=\Om$
and $\om_2=\pi^*\om_{11}+(A^\flat)^*(d\la)+d\la^{\om_{22},\na}$ and apply
Darboux's theorem. The key point is to check that $\om_1$ and $\om_2$ agree on the zero
section. This is a consequence of the expression of $(A^\flat)^*(d\la)$ computed in the
proof of Proposition \ref{propisition2} and the following formulas:
\begin{eqnarray*}
 d\la^{\om_{22},\na}(X^h,Y^h)(u)&=&\frac12\om_{22}(R(X,Y)u,u),\\
d\la^{\om_{22},\na}(X^v,Y^v)&=&\om_{22}(X,Y)\circ\pi,\\
d\la^{\om_{22},\na}(X^h,Y^v)(u)&=&\frac12\na_X\om_{22}(u,Y).
\end{eqnarray*}





\noindent Cadi-Ayyad University\\
Faculty of Science and Technology\\
BP 549 Marrakesh Morocco\\
Email:\\
abouqateb@fstg-marrakech.ac.ma\\
boucetta@fstg-marrakech.ac.ma\\
ikemakhen@fstg-marrakech.ac.ma\\

\end{document}